\documentstyle[12pt,psfig]{amsart}

\newtheorem{theorem}{Theorem}[section]

\newtheorem{algorithm}{Algorithm}[section] 

\newcommand{\nn}{{\Bbb N}}

\newcommand{\zz}{{\Bbb Z}}

\newcommand{\cc}{{\Bbb C}}
\newcommand{\pp}{{\Bbb P}}                 

\begin{document}

\bibliographystyle{plain}

\title[Polynomial Homotopies for dense, sparse and determinantal Systems]
      {Polynomial Homotopies for dense, sparse \\ and determinantal Systems}

\author{Jan Verschelde}
\address{\hskip-\parindent
         Mathematical Sciences Research Institute \\
         1000 Centennial Drive \\
         Berkeley, CA 94720-5070, U.S.A.} 
\address{\hskip-\parindent
         Department of Mathematics \\
         Michigan State University \\
         East Lansing, MI 48824-1027, U.S.A.} 

\email{jan@@msri.org, jan@@math.msu.edu or jan.verschelde@@na-net.ornl.gov,
and urls: http://www.msri.org/people/members/jan
http://www.mth.msu.edu/$\sim$jan}   

\thanks{Research at MSRI is supported in part by NSF grant DMS-9701755,
benefited from a post-doctoral fellowship at MSRI and is also
supported in part by NSF grant DMS-9804846 at MSU}

\date{July 11, 1999}

\begin{abstract}
\noindent 
Numerical homotopy continuation methods for three classes of polynomial
systems are presented.  For a generic instance of the class, every path
leads to a solution and the homotopy is optimal.
The counting of the roots mirrors the resolution of a generic system
that is used to start up the deformations.
Software and applications are discussed.

\bigskip


\bigskip

\end{abstract}

\maketitle

\tableofcontents 

\section{Introduction}

Solving polynomial systems numerically means computing approximations
to all isolated solutions.
Homotopy continuation methods provide paths to approximate solutions.
The idea is to break up the original system into simpler problems.
To solve the original system, the solutions of the simpler systems 
are deformed into the solutions of the original problem.

\medskip

This paper presents optimal homotopies for three different classes of 
polynomial systems.  Optimal means that for generic instances of the
classes there are no diverging solution paths, whence the amount of
computational work is linear in the number of solutions.
In the next section we list the principal key words, definitions and
main theorems for dense, sparse and determinantal polynomial systems.
The proofs of these theorems follow from the correctness of the homotopies.

\medskip

Path-following methods are standard numerical 
techniques~(\cite{ag90,ag93,ag97}, \cite{mor87}, \cite{wat86,wat89}) 
to achieve global convergence when solving nonlinear systems.
For polynomial systems we can reach all isolated solutions.
In the third section we describe
the paradigm of Cheater's homotopy~(\cite{lsy89},~\cite{lw92})
or coefficient-parameter polynomial continuation~(\cite{ms89},~\cite{ms90}).
This paradigm allows to construct homotopies for which singularities only
occur at the end of the paths. 
To deal with components of solutions we use an embedding
method that leads to generic points on each component.
This method is essential to numerical algebraic geometry~\cite{sw96}.

\medskip

From~\cite{kem93} we cite: ``Algebraic geometry studies the delicate balance
between the geometrically plausible and the algebraically possible''.
By a choice of coordinates we set up an algebraic formulation for a
geometric problem that is then solved by automatic computations.
While this approach is extremely powerful, we might get trapped into
tedious wasted computations after loosing the original geometric meaning
of the problem.  In section four we stress the geometric intuition of
homotopy methods.  Compactifications and homogeneous coordinates provide
us the tools to generate the numerically most favorable representations
for the solutions to our problem.
In section five we arrive at the heart of modern homotopy methods
where we outline specific algorithms to implement the root
counts\footnote{The term ``root count'' was coined by
Canny and Rojas~\cite{cr91} while introducing mixed volumes
to computational algebraic geometry.}.
The counting of the roots mirrors the resolution of a system in
generic position that is used as starting point in the deformations.

\medskip

Polyhedral methods occupy the central part of current research, as
they are responsible for a computational breakthrough in numerical
general-purpose solvers for polynomial systems.
Section six is devoted to numerical software with an emphasis on the
structure of the package PHC, developed by the author during the past decade.
Another novel and exciting research development concerns the numerical
Schubert calculus, which is one of the major new features in the second
public release of~PHC.
The author has gathered more than one hundred polynomial systems that
arose in various application fields.
This collection serves as a test suite for software and a gallery to
demonstrate the importance of polynomial systems to mathematical modelling. 
In section seven we sample some interesting cases.

\medskip

The reference list contains a compilation of the most relevant technical
contributions to polynomial homotopy continuation.  Besides those we
want to point at some other works in the literature that are of special
interest.  Some user-friendly introductions to algebraic geometry appeared
in recent years: see~\cite{abh90}, \cite{fal90}, \cite{har95}, with
computational aspects in~\cite{cls97} and~\cite{cls98}.
As Newton polytopes have become extremely important,
we recommend~\cite{zie95} and the handbook chapters~\cite{gr97}.
See also~\cite{stu98} for the interplay between the combinatorics of polytopes 
and the (real) roots of polynomials.
A recent survey that also covers polyhedral homotopies along with other
polynomial continuation methods appeared in~\cite{li97}.

\section{Three Classes of Polynomial Systems}

The classification in Table~\ref{tabclasses} is inspired by~\cite{hss98}.
The dense class is closest to the common algebraic description, whereas 
the determinantal systems arise in enumerative geometry.

\begin{table}[ht]
\begin{center}
\begin{tabular}{|c||c|c|c|c|} \hline
  system & model & theory & \multicolumn{2}{c|}{space} \\ \hline  \hline
  dense & highest degrees & B\'ezout & $\pp^n$ & projective \\ \hline
  sparse & Newton polytopes & Bernshte\v{\i}n & $(\cc^*)^n$ & toric \\ \hline
  determinantal & localization posets & Schubert & $G_{mr}$ & Grassmannian
                                                                    \\ \hline
\end{tabular}
\smallskip
\caption{Key words of the three classes of polynomial systems.}
\label{tabclasses}
\end{center}
\end{table}

For the vector of unknowns ${\bf x} = (x_1,x_2,\ldots,x_n)$
and exponents ${\bf a} = (a_1,a_2,\ldots,a_n) \in \nn^n$, denote 
${\bf x}^{\bf a} = x_1^{a_1} x_2^{a_2} \cdots x_n^{a_n}$.
A polynomial system $P({\bf x}) = {\bf 0}$ is given by
$P = (p_1,p_2,\ldots,p_n)$, a tuple of polynomials $p_i \in \cc[{\bf x}]$,
$i = 1,2,\ldots,n$.

\smallskip

The complexity of a {\em dense polynomial} $p$ is 
measured by its degree~$d$:
\begin{equation}
  p({\bf x}) = \sum_{0 \leq a_1+a_2+\cdots+a_n \leq d}
              c_{\bf a} {\bf x}^{\bf a}, \quad \quad d = \deg(p),
\end{equation}
where at least one monomial of degree $d$ should have a nonzero coefficient.
The {\em total degree} $D$ of a {\em dense system} $P$ is
$D = \prod_{i=1}^n \deg(p_i)$.

\begin{theorem}(B\'ezout~{\rm \cite{cls97}})
The system $P({\bf x}) = {\bf 0}$ has
no more than $D$ isolated solutions, counted with multiplicities.
\end{theorem}

Consider for example
\begin{equation} \label{eqexsys}
  P(x_1,x_2) =
  \left\{
    \begin{array}{r}
       x_1^4 + x_1 x_2 + 1 = 0 \\
       x_1^3 x_2 + x_1 x_2^2 + 1 = 0
    \end{array}
  \right.
  \quad {\rm with~total~degree}~D = 4 \times 4 = 16.
\end{equation}
Although $D = 16$, this system has only eight solutions because of its
sparse structure.

\medskip

The {\em support} $A$ of a {\em sparse polynomial} $p$ collects all
exponents of those monomials whose coefficients are nonzero.
Since we allow negative exponents (${\bf a} \in \zz^n$),
we restrict ${\bf x} \in (\cc^*)^n$, $\cc^* = \cc \setminus \{ 0 \}$.
\begin{equation}
  p({\bf x}) = \sum_{{\bf a} \in A} c_{\bf a} {\bf x}^{\bf a}, \quad \quad
  \forall {\bf a} \in A: c_{\bf a} \not= 0, \quad
  A \subset \zz^n, \quad \#A < \infty.
\end{equation}
The {\em Newton polytope} $Q$ of~$p$ is the convex hull of the support~$A$
of~$p$.
We model the structure of a {\em sparse system} $P$ by a tuple of Newton
polytopes ${\cal Q} = (Q_1,Q_2,\ldots,Q_n)$, spanned by 
${\cal A} = (A_1,A_2,\ldots,A_n)$, the so-called {\em supports} of~$P$.

\smallskip

\begin{figure}[hbt]
\begin{center}
\centerline{\psfig{figure=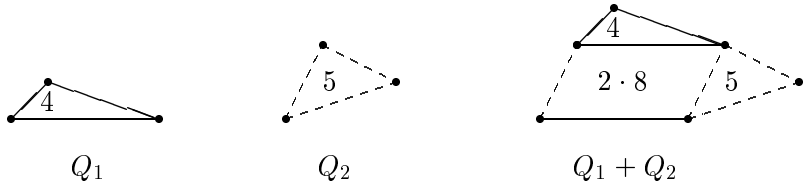}}
\caption{Newton polytopes $Q_1$, $Q_2$, a mixed subdivision
of $Q_1 + Q_2$ with volumes.}
\label{figmcc}
\end{center}
\end{figure}

The volume of a positive linear combination of polytopes is a homogeneous
polynomial in the multiplication factors.
The coefficients are {\em mixed volumes}.
For instance, for $(Q_1,Q_2)$, we write:
\begin{equation} \label{eqvolpol}
  2!{\rm vol}_2(\lambda_1 Q_1 + \lambda_2 Q_2)
   =  V_2(Q_1,Q_1) \lambda_1^2  + 2 \cdot V_2(Q_1,Q_2) \lambda_1 \lambda_2
      + V_2(Q_2,Q_2) \lambda_2^2,
\end{equation}
normalizing $V_2(Q,Q) = 2! {\rm vol}_2(Q)$.
For the Newton polytopes of the system~(\ref{eqexsys}):
$2! {\rm vol}_2(\lambda_1 Q_1 + \lambda_2 Q_2) =
 4 \lambda_1^2 + 2 \cdot 8 \lambda_1 \lambda_2 + 5 \lambda_2^2$.
To interpret this we look at Figure~\ref{figmcc} and see that multiplying
$P_1$ and $P_2$ respectively by $\lambda_1$ and $\lambda_2$ changes their
areas respectively with $\lambda_1^2$ and $\lambda_2^2$.
The cells in the subdivision of $Q_1 + Q_2$ whose area is scaled by
$\lambda_1 \lambda_2$ contribute to the mixed volume.
So, for the example in~(\ref{eqexsys}), the root count is eight.

\begin{theorem} (Bernshte\v{\i}n~{\rm \cite{ber75}})
A system $P({\bf x}) = {\bf 0}$ with Newton polytopes~${\cal Q}$ has no
more than $V_n({\cal Q})$ isolated solutions in $(\cc^*)^n$,
counted with multiplicities.
\end{theorem}
The mixed volume was nicknamed~\cite{cr91} as the BKK bound to honor
Bernshte\v{\i}n~\cite{ber75}, Kushnirenko~\cite{kus76}, and
Khovanski{\v{\i}}~\cite{kho78}.

\medskip

For the third class of polynomial systems we consider a matrix $[C | X]$
where $C \in \cc^{(m+r) \times m}$ and $X \in \cc^{(m+r) \times r}$
respectively collect the coefficients and indeterminates.
Laplace expansion of the maximal minors of $[C | X]$ in $m$-by-$m$ and
$r$-by-$r$ minors yields a {\em determinantal} polynomial 
\begin{equation} \label{eqdetpol}
  p({\bf x}) = \sum_{\scriptsize \begin{array}{c} I \cup J = U \\
									  I \cap J = \emptyset
                     \end{array} } {\rm sign}(I,J) C[I] X[J],
  \quad U = \{ 1,2,\ldots,m+r\},
\end{equation}
where the summation runs over all distinct choices $I$ of $m$ elements of~$U$.
The partition $\{ I, J \}$ of $U$ defines the permutation $U \mapsto (I,J)$
with ${\rm sign}(I,J)$ its sign.
The symbols $C[I]$ and $X[I]$ respectively represent
coefficient minors and minors of indeterminates.
Note that for more general intersection conditions, the matrices
$[C | X]$ are not necessarily square.

\smallskip

The vanishing of a polynomial as in~(\ref{eqdetpol}) expresses
the condition that the $r$-plane $X$ meets a given $m$-plane nontrivially.
The counting and finding of all figures that satisfy certain geometric
conditions is the central theme of enumerative geometry.
For example, consider the following.

\begin{theorem} (Schubert~{\rm \cite{sch1891}})
Let $m,r \geq 2$.  In~$\cc^{m+r}$ there are
\begin{equation} \label{eqgrassdeg}
   d_{m,r} \quad = \quad
   \frac{1! \, 2! \, 3! \cdots (r\!- \!2) ! \, (r \!-\!1)! \cdot
     (mr)!}{m!\, (m \! + \! 1)! \, (m \! + \! 2)!
     \cdots(m \! + \! r \! - \! 1)!}
\end{equation}        
$r$-planes that nontrivially meet $mr$ given $m$-planes in general position.
\end{theorem}
This root count $d_{m,r}$ is sharp compared to other root counts,
see~\cite{sot98} and~\cite{ver98} for examples.

\smallskip

We can picture the simplest case, using the fact that 2-planes in $\cc^4$
represent lines in $\pp^3$.  In Figure~2 the positive real projective 3-space
corresponds to the interior of the tetrahedron.

\bigskip

\parbox[b]{5cm}{

\centerline{\psfig{figure=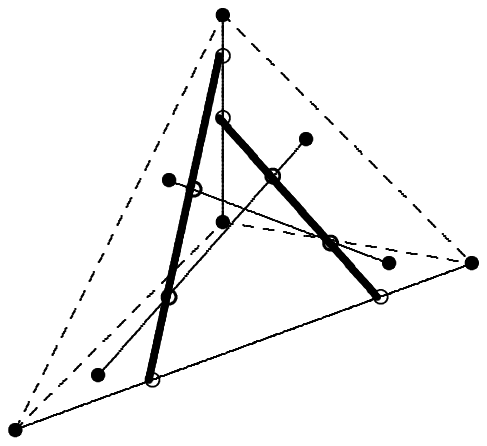,width=5cm}}

~~~~~Figure~2: $m = 2 = r$.

\addtocounter{figure}{1}

} \ \ \ \ \ \ \ \
\parbox[b]{10cm}{ 

\smallskip

In Figure~2, we see two thick lines meeting four given skew fine lines
in a point.  When not all input planes have the same dimension,
but when the number of solutions is still finite,
Pieri's formula~\cite{pie1891} provides a root count~\cite{sot97a,sot97b}.

\smallskip

In~\cite{hss98} the problem is solved in chains of nested subspaces,
using a cellation of the Grassmannian $G_{mr}$ of $r$-planes in $\cc^{m+r}$.
A {\em localization poset} models~\cite{hv99} the specialization of
the solution $r$-plane when the input is specialized.

}

\bigskip

Algorithmic proofs for the above theorems consist in two steps.
First we show how to construct a generic start system that has exactly
as many regular solutions as the root count.
Then we set up a homotopy for which all isolated solutions of any
particular target system lie at the end of some solution path 
originating at some solution of the constructed start system.

\section{The Principles of Polynomial Homotopy Continuation Methods}

Homotopy continuation methods operate in two stages.
Firstly, homotopy methods exploit the structure of $P$ to find a root count 
and to construct a start system $P^{(0)}({\bf x}) = {\bf 0}$
that has exactly as many regular solutions as the root count. 
This start system is embedded in the {\em homotopy}
\begin{equation} \label{eqlinhom}
   H({\bf x},t) = \gamma(1-t)P^{(0)}({\bf x}) + t P({\bf x}) = {\bf 0},
   \quad t \in [0,1],
\end{equation}
with $\gamma \in \cc$ a random number.
In the second stage, as $t$ moves from 0 to 1,
numerical continuation methods trace the paths that originate at the 
solutions of the start system towards the solutions of the target system.

\smallskip

\noindent The good properties we expect from a homotopy
$H({\bf x},t) = {\bf 0}$ are (borrowed from \cite{li97}):
\begin{enumerate}
 \item ({\em triviality}) The solutions for~$t=0$ are trivial to find.
 \item ({\em smoothness}) No singularities along the solution paths occur.
 \item ({\em accessibility}) All isolated solutions can be reached.
\end{enumerate}

\bigskip

Continuation or path-following methods are standard numerical techniques
(\cite{ag90,ag93,ag97}, \cite{mor87}, \cite{wat86,wat89}) to trace the
solution paths defined by the homotopy using {\em predictor-corrector} methods.
The smoothness property of complex polynomial homotopies implies that paths
never turn back, so that during correction the parameter $t$ stays fixed,
which simplifies the set up of path trackers. 
A pseudo-code description of a path tracker is in
Algorithm~\ref{algpathfoll}.

\medskip

\begin{figure}[hbt]
\begin{center}
\centerline{\psfig{figure=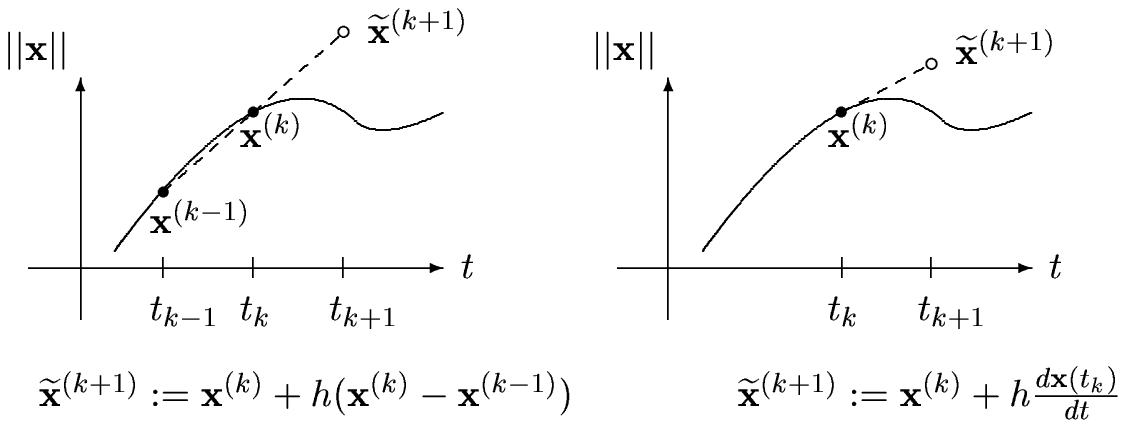}}
\caption{The secant and tangent predictor with step length~$h$.}
\label{figsectan}
\end{center}
\end{figure}

The {\em predictor} delivers at each step of the method
a new value for the continuation parameter and predicts an approximate
solution of the corresponding new system in the homotopy.
Figure~\ref{figsectan} shows two common predictor schemes.
The predicted approximate solution is adjusted by applying Newton's method
as {\em corrector}.  The third ingredient in path-following methods is
the {\em adaptive step size control}.
The step length is determined to enforce quadratic convergence in the
corrector to avoid path crossing.        

\medskip

\begin{algorithm} \label{algpathfoll}
{\rm Following one solution path by an increment-and-fix
predictor-corrector method with an adaptive step size control strategy.

\bigskip

\noindent \begin{tabular}{lcr}
  Input: \ \ $H({\bf x},t)$, ${\bf x}^* \in \cc^n$: $H({\bf x}^*,0) = {\bf 0}$,
               & \ \ \ & {\em homotopy and start solution} \\
  \ \ \ \ \ \ \ \ \ \ \ \ $\epsilon > 0$, $max\_it$, $max\_steps$.
               & & {\em accuracy and upper bounds} \\
  Output:  ${\bf x}^*$, success if $||H({\bf x}^*,1)|| \leq \epsilon$.
               & & {\em approximate solution if success} \\
\\
  $t := 0$; \ $k := 0$;     & & {\em initialization} \\
  $h := max\_step\_size$;   & & {\em step length} \\
  $old\_t := t$; \ $old\_{\bf x}^* := {\bf x}^*$
                            & & {\em back up values for $t$ and ${\bf x}^*$ } \\
  $previous\_{\bf x}^* := {\bf x}^*$; 
                            & & {\em previous approximate solution} \\
  stop := false;            & & {\em combines stopping criteria} \\ 
  while $t < 1$ and not stop loop & & \\
\ \ \ $t := \min(1,t + h)$; & & {\em secant predictor for $t$} \\
\ \ \ ${\bf x}^* := {\bf x}^* + h ( {\bf x}^* - previous\_{\bf x}^* )$;
                            & & {\em secant predictor for ${\bf x}^*$} \\
\ \ \ Newton($H({\bf x},t),{\bf x}^*,\epsilon,max\_it$,success);
                            & & {\em correct with Newton's method} \\
\ \ \ if success            & & {\em step size control} \\
\ \ \ \ then $h := \min(Expand(h),max\_step\_size)$; 
                            & & {\em enlarge step length} \\
\ \ \ \ \ \ \ \ \ \ \ $previous\_{\bf x}^* := old\_{\bf x}^*$; 
                            & & {\em go further along path} \\
\ \ \ \ \ \ \ \ \ \ \ $old\_t := t$; \ $old\_{\bf x}^* := {\bf x}^*$;
                            & & {\em new back up values} \\
\ \ \ \ else \ $h := Shrink(h)$; & & {\em reduce step length} \\
\ \ \ \ \ \ \ \ \ \ \ $t := old\_t$; \ ${\bf x}^* := old\_{\bf x}^*$;
                            & & {\em step back and try again} \\
\ \ \ end if; & & \\
\ \ \ $k := k+1$;           & & {\em augment counter} \\
\ \ \ stop := ($h < min\_step\_size$) or ($k > max\_steps$);
                            & & {\em stopping criteria} \\
  end loop; & & \\
  success := ($||H({\bf x}^*,1)|| \leq \epsilon$).
                            & & {\em report success or failure}
\end{tabular}
}
\end{algorithm}

Following all paths can be done sequentially, one path at a time, or
in parallel, with for each solution path the same sequence of values of
the continuation parameter.
The sequential path-following method has the advantage that the low
overhead of communication~\cite{acw89} makes it very suitable to run on
multi-processor environments.
Note that the memory requirements are optimal. 

\medskip

To solve repeatedly a polynomial system with the same coefficient structure
$P({\bf c},{\bf x}) = {\bf 0}$, the homotopy~(\ref{eqlinhom}) is applied
with $P^{(0)} = P({\bf c}^0,{\bf x}) = {\bf 0}$ a system with random
coefficients~${\bf c}^0$.
Solving $P({\bf c}^0,{\bf x}) = {\bf 0}$ is no longer trivial,
so the name {\em cheater's homotopy}~\cite{lsy89} is appropriate.
A similar idea appeared in~\cite{ms89,ms90}.
For coefficients given as functions of parameters, a refined version of 
cheater's homotopy in~\cite{lw92} avoids repeated evaluation of those
functions during path following:
\begin{equation} \label{eqnonlincheat}
  H({\bf x},t) = P((1-[t-t(1-t) \gamma]){\bf c}^0
                 + (t-t(1-t) \gamma){\bf c},{\bf x}) = {\bf 0},
                 \quad t \in [0,1], \gamma \in \cc.
\end{equation}
In~\cite{lw92} it is proven that with~(\ref{eqnonlincheat}) all isolated
solutions of~$P({\bf c},{\bf x}) = {\bf 0}$ can be reached and that
singularities can only occur at the end of the paths.

\smallskip

Typically, when using a cheater's homotopy, the computational effort spent
towards the end of the paths often accounts for most of the work.
The main numerical problem is then to distinguish irrelevant solutions
at infinity from ill-conditioned but possibly meaningful solutions.
End games~\cite{hv98}, \cite{msw91,msw92a,msw92b}, \cite{sws96}
provide several procedures to approximate the winding number of a path.
Recently, Zeuthen's rule was applied in~\cite{kss98} to determine
numerically the multiplicity of an isolated solution.
Multi-precision facilities are useful for evaluation of residuals and
root refinement for badly scaled solutions. 

\medskip

In most applications, the polynomial systems have real coefficients and
invite the use of real homotopies.
In~\cite{bm84} it was conjectured and proven in~\cite{lw93}
that generically, real homotopies contain no singular points other
than a finite number of quadratic turning points.
At those bifurcation points pairs of real solution paths become imaginary
or conversely, complex conjugated solution paths join to yield two real
solution paths.
We refer to~\cite{all84}, \cite{hk90}, \cite{li97} and~\cite{lw93,lw94}
for a discussion of numerical techniques to deal with quadratic turning points.
A remarkable application of real homotopies in the real world
consists in the finding of the relevant parameters of a polynomial system
to maximize the number of real roots, see~\cite{die98}
for the 40 real solutions for the Stewart-Gough platform in mechanics.

\medskip

In~\cite{sw96} the use of homotopy continuation to deal with overdetermined
and components of solutions is discussed.
Geometrically one slices the components of solutions with as many random
hyperplanes as the dimension of the components.
The solutions to the original polynomial system augmented with these
random linear equations for the hyperplanes are {\em generic points}
of the components, constituting the main numerical data to study those
components.  In particular, the number of generic points one obtains
by this slicing procedure equals the sum of the degrees over all
top-dimensional components of solutions.

\medskip

To make the algorithms of~\cite{sw96} more efficient, in~\cite{sv99},
the following embedding of the polynomial system~$P({\bf x}) = {\bf 0}$
is proposed:
\begin{equation} \label{eqembed}
  \left\{
     \begin{array}{rl}
        p_i({\bf x}) + \lambda_i z = 0, & i=1,2,\ldots,n \\
        {\displaystyle \sum_{j=1}^n c_j x_j + z = 0}
     \end{array}
  \right.
\end{equation}
where the $\lambda_i$'s and $c_j$'s are random complex numbers.
This embedding has the advantage over the algorithms in~\cite{sw96} that
fewer solution paths diverge.  Solutions to the system~(\ref{eqembed})
with $z = 0$ lie on a component of solutions.
By Bertini's theorem, all solutions with $z \not= 0$ are regular.
In~\cite{sv99}, it is proven that those solutions can be used as
start solutions to reach {\em all} isolated solutions of the original
polynomial system~$P({\bf x}) = {\bf 0}$.

\medskip

The embedding~(\ref{eqembed}) is performed repeatedly in the routine
`Embed' in the algorithm (copied from~\cite{sv99}) below.

\begin{algorithm} \label{algcascade}
{\rm Cascade of homotopies between embedded systems.

\bigskip

\begin{tabular}{lcr}
Input: $P$, $n$.
    & & {\em system with  solutions in $\cc^n$} \\
Output: $({\cal E}_i,{\cal X}_i,{\cal Z}_i)_{i=0}^n$.
    & & {\em embeddings with solutions} \\
\\
${\cal E}_0 := P$;
    & & {\em initialize embedding sequence} \\
for $i$ from 1 up to $n$ do
    & & {\em slice and embed} \\
\ \ \ ${\cal E}_i$ := Embed(${\cal E}_{i-1},z_{i}$);
    & & {\em $z_i$ = new added variable} \\
end for; & & {\em homotopy sequence starts} \\
${\cal Z}_n$ := Solve(${\cal E}_n$);
    & & {\em all roots are isolated, nonsingular, with $z_n \not= 0$}\\
for $i$ from $n-1$ down to 0 do 
    & & {\em countdown of dimensions} \\
\ \ \ $H_{i+1}$ := $
  t {\cal E}_{i+1}
  +
  (1-t)\left(
     \begin{array}{c}
        {\cal E}_i \\
        z_{i+1}
     \end{array}
  \right) $;
& &  {\em \begin{tabular}{r}
              homotopy continuation \\
       $t: 1 \rightarrow 0$ to remove $z_{i+1}$ \\
         \end{tabular} } \hspace{-5mm} \\
 \ \ \ ${\cal X}_i$ := limits of solutions of $H_{i+1}$ \\
 \ \ \ \ \ \ as $t\to 0$ with $z_i=0$;
    & &  {\em on component} \\
\ \ \ ${\cal Z}_i$ := $H_{i+1}({\bf x},z_{i} \not= 0,t=0)$;
    & &  {\em not on component: these solutions} \\
    & &  {\em  are  isolated and nonsingular} \\
end for. & & \\
\end{tabular}
}
\end{algorithm}

This embedding allows the efficient treatment of overdetermined systems
and other nonproper intersections.
By perturbing the added hyperplanes and extending the generic points by
continuation, interpolation methods can lead to equations for the components.

\section{The Geometry of the Deformations}

Homotopy methods have an intuitive geometric interpretation.
In this section we illustrate the geometry of the three types of
moving into special position: product, toric, and Pieri deformations.
These can be regarded as three applications of the principle of 
continuity or conservation of number in enumerative geometry.

\smallskip

Product homotopies deform polynomial equations into products of
linear equations.  In Figure~\ref{figbasic} we see the line configuration
at the start and the ellipse-parabola intersection in the end.
Note that complex space is the natural space for deformations.
The other two complex conjugated intersection points could not be
displayed in Figure~\ref{figbasic}.

{\small
\begin{figure}[hbt]
\begin{center}
\centerline{\psfig{figure=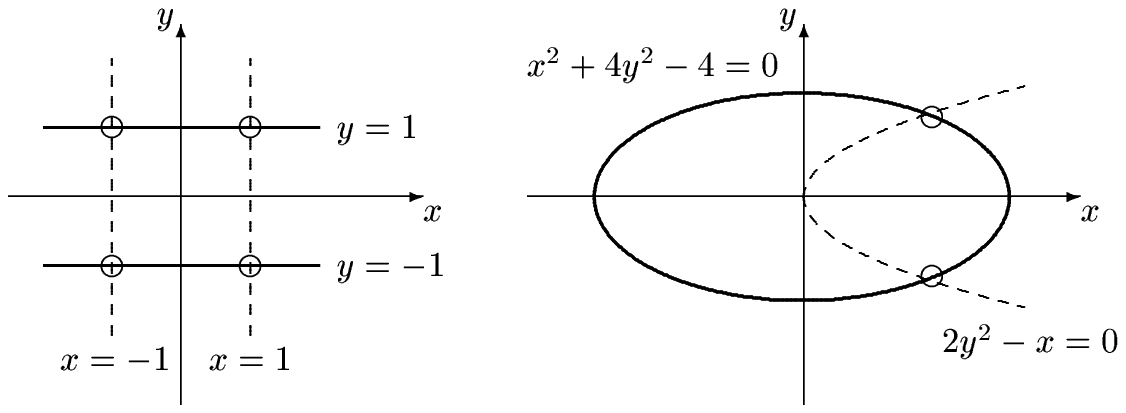}}
\caption{Intersection of quadrics: a degenerate and a target configuration.}
\label{figbasic}
\end{center}
\end{figure}
}

\medskip

The sparser a system, the easier it can be solved.
In Figure~\ref{figincrpoco} we illustrate the idea of making a system
sparser by setting up a so-called polyhedral homotopy that reduces this
particular system at $t=0$ to a linear system.
The lower hull of the Newton polytope of this homotopy
induces a triangulation, which is used to count the roots.
In particular, every cell in the triangulation gives rise to a homotopy
with as many paths to follow as the volume of the cell.
The other root for the example in Figure~\ref{figincrpoco} can be computed
with a homotopy obtained from~$\widehat P$ by the substitution 
of variables $x_1 \leftarrow {\widetilde x}_1 t^{-1}$ and
$x_2 \leftarrow {\widetilde x}_2 t^{-1}$.  This transformation pushes the
constant monomial up, so that at $t=0$ we have the nonconstant monomials
in the start system to compute the other root.

{\small
\begin{figure}[hbt]
\begin{center}
\centerline{\psfig{figure=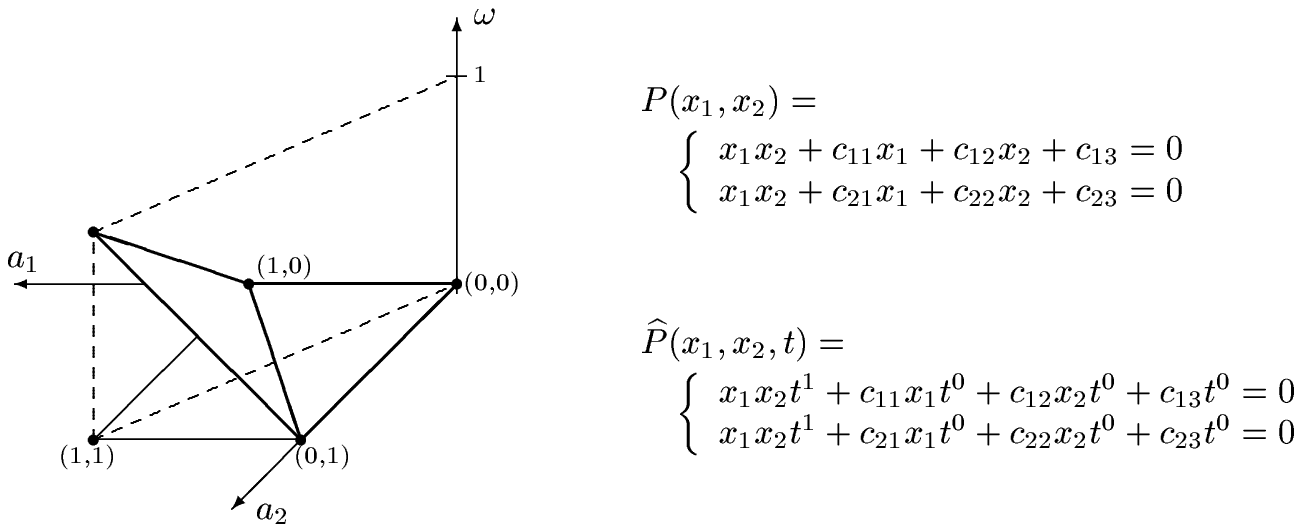}}
\caption{Triangulation of the Newton polytope of $P$ with
polyhedral homotopy $\widehat P$.}
\label{figincrpoco}
\end{center}
\end{figure}
}

\medskip

\begin{figure}[hbt]
\begin{center}
\centerline{\psfig{figure=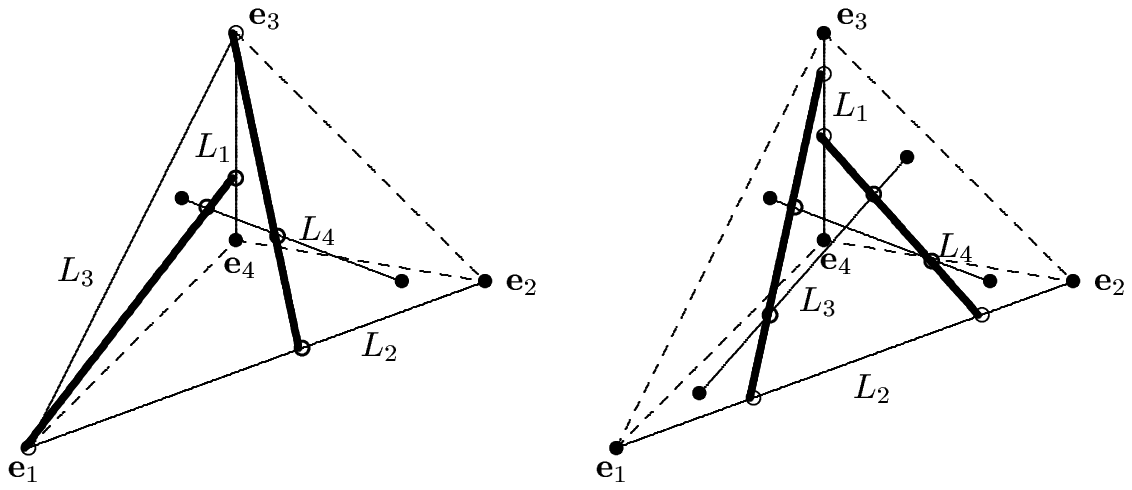}} 
\caption{In $\pp^3$ two thick lines meet four given lines
$L_1$, $L_2$, $L_3$, and $L_4$ in a point.
At the left we see a special configuration and the general
configuration is at the right.}
\label{figlines}
\end{center}
\end{figure}

Figure~\ref{figlines} displays a special and a general configuration
of four lines.  The basis has been chosen such that two of the four
input lines are spanned by standard basis vectors.
To compute all lines that meet four given lines, one of the four
given lines is moved into special position so that it intersects
two other given lines, see the left of Figure~\ref{figlines}.
The solution lines must then originate at those two intersection points
and reach to the other opposite line while meeting the line left in
general position.

\bigskip

The constructions above are in a sense~\cite{abh90} ``heuristic proofs''.
With the general position assumption we cheat a bit, avoiding the hard
problem of assigning multiplicities.  Making this so-called~\cite{wei62}
``method of degeneration'' rigorous was an important development in 
algebraic geometry.

\bigskip

To deal with solution paths diverging to ill-conditioned
roots or to infinity we need to compactify our space.
Instead of polynomials in $n$ variables we consider homogeneous forms
with coordinates subject to equivalence relations.
While mathematically all coordinate choices are equivalent, 
we select the numerically most favorable representations
of the solutions.

\medskip

The usual projective transformation consists in the change of
variables $x_i := \frac{z_i}{z_0}$, for $i=1,2,\ldots,n$, which
leads to the homogeneous system $P({\bf z}) = {\bf 0}$.
To have as many equations as unknowns, we add to this system a
random hyperplane.  Except for an algebraic set of the coefficients
of this added hyperplane, all solution paths are guaranteed to stay
inside $\cc^{n+1}$ when homotopy continuation is applied.
We refer to~\cite{li97} for numerical techniques that dynamically
restrict the computations to $n$ dimensions.

\medskip

For sparse polynomial systems, we introduce as in~\cite{ver99}
a new variable for every facet of the Newton polytopes.
The advantage of this more involved compactification
is based on the observation that when paths diverge
to infinity certain coefficients of the polynomial system become dominant.
With toric homogenization the added variables that become zero identify
the faces of the Newton polytopes for the parts of the system that
become dominant.  This compactification works in conjunction with 
polyhedral end games~\cite{hv98} which are summarized in 
Section~\ref{secsparhom}.

\medskip

The natural way to compactify $G_{mr}$ is to consider a multi-projective
homogenization according to rows or columns of the matrix representations
for the planes.  In addition, we have that the planes are equivalent upon
a linear change of basis.  Choosing orthonormal matrices to represent the
input planes leads to drastic improvements in the conditioning of the
solution paths, see~\cite{hv99} and~\cite{ver98} for experimental data.

\section{Root Counts and Start Systems}

The main principle is that counting roots corresponds to solving start
systems.  Algorithms to illustrate this principle will be shown for little
examples for the three classes of polynomial systems.

\smallskip

For dense polynomial systems, the computation of generalized permanents
model the resolution of linear-product start systems.
The algorithms to compute mixed volumes lead to polyhedral homotopies
to solve sparse polynomial systems.
The localization posets describe the structure of the cellation of
the Grassmannian used to set up the Pieri deformations.

\subsection{Dense Polynomials modeled by Highest Degrees}

A polynomial in one variable has as many complex solutions as its degree.
A linear system has either infinitely many solutions or 
exactly one isolated solution in projective space.
By this analogy~\cite{li87} we see that B\'ezout's theorem generalizes
these last two statements:
a polynomial system has either infinitely many solutions or exactly
as many isolated solutions in complex projective space as the total degree.

\medskip

As the above presentation of B\'ezout's theorem suggests, the simplest
cases are univariate and linear systems, which are used as start systems.
For the example system~(\ref{eqexsys}),
a start system~$P^{(0)}({\bf x}) = {\bf 0}$ based on the
total degree~$D$ is given by two univariate quartic equations
$x_1^4 - c_1 = 0 = x_2^4 - c_2$,
where $c_1$ and $c_2$ are randomly chosen complex numbers.
Note that the computation of $D = 4 \times 4$ models the structure of the
solutions of $P^{(0)}$ as four solutions for $x_1$ crossed with four solutions 
for~$x_2$.

\medskip

The earliest approaches of this homotopy
appear in~\cite{cmy79}, \cite{dre77}, \cite{gz79}, \cite{gl80},
and were further developed in~\cite{li83}, \cite{mor83}, \cite{wri85}.
The book~\cite{mor87} contains a very good introduction to the
practice of solving polynomial system by homotopy continuation.
Regularity results can be found in~\cite{ls87} and~\cite{zul88}.
While this homotopy algorithm has a sound theoretical basis,
the total degree is a too crude upper bound for the number of 
affine roots to be useful in most applications.

\medskip

Multi-homogeneous homotopies were introduced in~\cite{ms87a,ms87b} and applied
to various problems in mechanism design, see e.g.~\cite{wms90,wms92a}.
Similar are the random product homotopies~\cite{lsy87b,lsy87a},
applying intersection theory in~\cite{lw91},
but less suitable for automatic procedures.
For our running example~(\ref{eqexsys}),
we follow the approach of multi-homogenization
and we group the unknowns according to the partition
$Z = \{ \{ x_1 \} , \{ x_2 \} \}$.  The corresponding degree matrix $M_Z$
has in its $(i,j)$-th entry the degree of the $i$-th polynomial in the
variables of the $j$-th set of $Z$.
The 2-homogeneous B\'ezout bound $B_Z$ is the permanent of~$M_Z$.
\begin{equation}
   \begin{array}{l}
      P(x_1,x_2) = \\
      ~~~~~ \left\{
        \begin{array}{r}
           x_1^4 + x_1 x_2 + 1 = 0 \\
           x_1^3 x_2 + x_1 x_2^2 + 1 = 0
        \end{array}
      \right.
   \end{array}
\quad
  \begin{array}{rcc}
   M_Z & \! \! = \! \! & M_{\{ \{ x_1 \} , \{ x_2 \} \}} \\
       & \! \! = \! \! & \left[
       \begin{array}{cc}
          4 & 1 \\
          3 & 2
       \end{array}
    \right]
  \end{array}
\quad
    \begin{array}{rcl}
      B_Z & \! \! = \! \! & {\rm per}(M_Z) \\
          & \! \! = \! \! & 4 \times 2 + 3 \times 1 \\
          & \! \! = \! \! & 11
    \end{array}
\end{equation}
The computation of the permanent follows the expansion for the determinant,
except for the permanence of signs, as it corresponds to adding up the roots
when solving the corresponding linear-product start system:
\begin{equation}
  P^{(0)}({\bf x}) =
    \left\{
      \begin{array}{c}
         {\displaystyle \prod_{i=1}^4 (x_1 - \alpha_{1i}) 
                        \prod_{i=1}^1 (x_2 - \beta_{1i})} = 0 \\
         {\displaystyle \prod_{i=1}^3 (x_1 - \alpha_{2i})
                        \prod_{i=1}^2 (x_2 - \beta_{2i})} = 0 \\
      \end{array}
    \right.
\end{equation}
In most applications the grouping of variables follows from their meaning,
e.g.: for eigenvalue problems~$A {\bf x} = \lambda {\bf x}$,
 $Z = \{ \{ \lambda \} , \{ x_1, x_2 ,\ldots ,x_n \} \}$.
Efficient permanent evaluations in conjunction with
exhaustive searching algorithms for finding an optimal grouping
were developed in~\cite{wam92}.  In case the number of independent roots
equals the B\'ezout bound, interpolation methods~\cite{vbh91} are useful.

\medskip

Partitioned linear-product start systems were developed in~\cite{vh93a}
elaborating the idea that several different partitions can be used for
the polynomials in the system.
Motivated by symmetric applications~\cite{vc94a}, general linear-product
start system were proposed in~\cite{vc93b}.  These start systems are
based on a supporting set structure $S$ which provides a more refined model
of the degree structure of a polynomial system.
\begin{equation}
 S = 
  \begin{array}{|c|} \hline
     \{ x_1 \} \{ x_1 , x_2 \} \{ x_1 \} \{ x_1 \} \\ 
     \{ x_1 \} \{ x_1 , x_2 \} \{ x_1 \} \{ x_2 \} \\ \hline
  \end{array}
\end{equation}
To compute the bound formally, one collects all admissible $n$-tuples of sets,
picking one set out of every row in the set structure.
\begin{equation}
  \begin{array}{rl}
     B_S = & \! \! \! \! \# \{ ( \{ x_1 \} , \{ x_1 , x_2 \} ) ,
                   ( \{ x_1 \} , \{ x_2 \} ) ,
                   ( \{ x_1 , x_2 \} , \{ x_1 \} ) , \\
           &  ~~~  ( \{ x_1 , x_2 \} , \{ x_1, x_2 \} ) , 
                   ( \{ x_1 , x_2 \} , \{ x_1 \} ) , 
                   ( \{ x_1 , x_2 \} , \{ x_2 \} ) , \\
           &  ~~~  ( \{ x_1 \} , \{ x_2 \} ) , 
                   ( \{ x_1 , x_2 \} , \{ x_2 \} )
                   ( \{ x_1 \} , \{ x_2 \} )
                   ( \{ x_1 \} , \{ x_2 \} ) \}
  \end{array}
\end{equation}
Each admissible pair corresponds to a linear system that leads to a
solution of a generic start system:
\begin{equation}
  P^{(0)}({\bf x}) =
  \left\{
    \begin{array}{c}
       (x_1+c_{11})(x_1+c_{12} x_2+c_{13})(x_1+c_{14})(x_1+c_{15}) = 0 \\
       (x_1+c_{21})(x_1+c_{22} x_2+c_{23})(x_1+c_{24})(x_2+c_{25}) = 0 \\
    \end{array}
  \right.    
\end{equation}
This start system has $B_S = 10$ solutions.
In~\cite{vc93b}, the following theorems were proven.

\begin{theorem} Except for a choice of coefficients belonging to an
algebraic set, there are exactly $B_S$ regular solutions to a random
linear-product system based on the set structure~$S$.
\end{theorem}
The proof of the theorem consists in collecting the determinants
that express the degeneracy conditions.  These determinants are
polynomials in the coefficients and vanish at an algebraic set.

\medskip

\begin{theorem} All isolated solutions to~$P({\bf x}) = {\bf 0}$ lie
at the end of some solution path defined by a convex-linear homotopy
originating at a solution of a random linear-product start system,
based on a supporting set structure for~$P$.
\end{theorem}
The idea of the proof is to embed the homotopy into an appropriate
projective space and to consider the projection of the discriminant
variety as an algebraic set for the continuation parameter.
See~\cite{lww96} for an alternative proof.

\medskip

A general approach to exploit product structures was developed
in~\cite{msw95}.  For systems whose polynomials are sums of products
one may arrive at a much tighter bound replacing the products by one
simple product.  An efficient homotopy to solve the nine-point problem
in mechanical design was obtained in this way.

\medskip

The complexity of this homotopy based on the total degree is 
addressed in~\cite{bcss97} where $\alpha$-theory is applied
to give bounds on the number of steps that is needed to trace
the solution paths.  A major result is that one can decide in
polynomial time whether an average polynomial system has a solution.
A similar analysis of Newton's method in multi-projective space
was recently done in~\cite{ds97}.

\medskip

While the above complexity results apply to random systems,
the problem of automatically extracting and exploiting the
degree structure of a polynomial system is a much harder problem.
Finding an optimal multi-homogeneous grouping essentially requires
the enumeration of all partitions~\cite{wam92}.
With supporting set structures one may obtain a high success rate,
see~\cite{lww96} for a efficient heuristic algorithm.
Recent software extensions for finding optimal partitioned linear-product
start systems are in~\cite{wsw98}.

\subsection{Mixed Subdivisions of Newton Polytopes to compute Mixed Volumes}

For~(\ref{eqexsys}), we collect the exponent vectors of the system
$P$ in the supports $\cal A$:

\begin{equation} \label{eqexpolysup}
  \begin{array}{l}
   P(x_1,x_2) = \\
   ~~~~~ \left\{
       \begin{array}{r}
          x_1^4 + x_1 x_2 + 1 = 0 \\
          x_1^3 x_2 + x_1 x_2^2 + 1 = 0
       \end{array}
     \right.
  \end{array}
   \quad 
  \begin{array}{l}
	 {\cal A} = (A_1,A_2) \\
     ~~~~ A_1 = \{ (0,0) , (1,1) , (4,0) \} \\
     ~~~~ A_2 = \{ (0,0) , (1,2) , (3,1) \}
  \end{array}
\end{equation}
The supports $A_1$ and $A_2$ span the respective Newton polytopes
$Q_1$ and $Q_2$.

\medskip

The Cayley trick~\cite[Proposition 1.7, page 274]{gkz94} is a method
to rewrite a certain resultant as a discriminant of one single
polynomial with additional variables.
The polyhedral version of this trick as in~\cite[Lemma 5.2]{stu94}
is due to Bernd Sturmfels. 
It provides a one-to-one correspondence between the cells in a mixed
subdivision and a triangulation of the so-called Cayley polytope spanned
by the points of $A_i$ embedded in a $(2n-1)$-dimensional space.
See~\cite{hrs99} for another application besides mixed-volume computation.
As in~\cite{hrs99}, Figure~7 gives a ``one-picture proof'' of this trick,
displaying the Cayley polytope for the supports~${\cal A}$
in~(\ref{eqexpolysup}).  Note that this construction provides a
definition for mixed subdivisions.

\medskip

\noindent~\parbox[t]{8cm}{

The Cayley polytope is spanned by the points in $A_1$,
where each point of $A_1$ is extended with $n-1$ zero coordinates,
and the points in $A_i$ where each point in $A_i$ is extended
with the respective $i$-th standard basis vector, for $i=1,2,\ldots,n-1$.

\smallskip

Omitting the added coordinates of this Cayley embedding, every cell
in a triangulation of the Cayley polytope is identified with a cell
in a mixed subdivision of the original tuple of polytopes.
We can see this identification geometrically when slicing the
Cayley polytope with a hyperplane that separates the embedded
polytopes.  As in Figure~7, the slice contains 
$\lambda_1 Q_1 + \lambda_2 Q_2$ and the cells of a mixed subdivision
are cut out by the cells in a triangulation of the Cayley polytope.

}~\hspace{-0.3cm} 
\parbox[t]{10cm}{
\begin{center}
\centerline{\psfig{figure=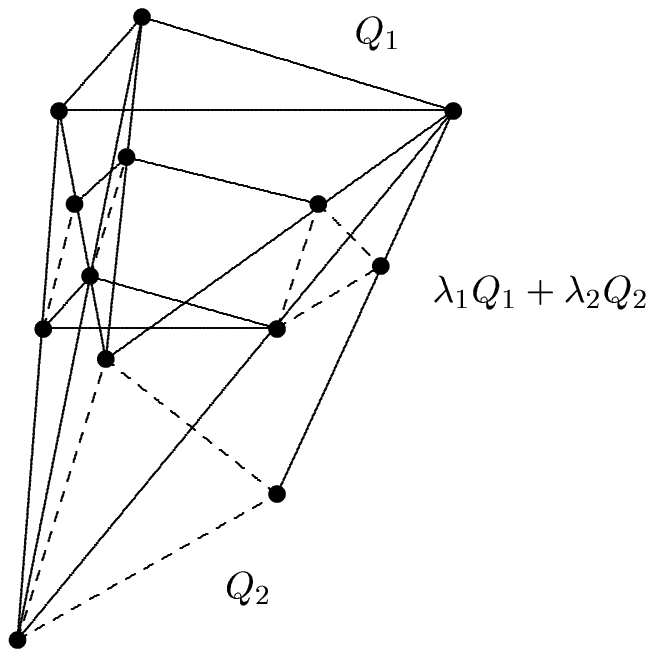}}

\medskip

Figure~7 : Cayley polytope of $Q_1$ and $Q_2$.
\addtocounter{figure}{1}
\end{center}

}

\medskip

The Cayley trick was
implemented in~\cite{vgc96} as an application of the dynamic
lifting algorithm to construct regular triangulations.
This method calculates the volume polynomial~(\ref{eqvolpol}) completely.
When one is only interested in the mixed volume, the method is only
efficient when the supports do not differ much from each other.

\medskip

\begin{figure}[hbt]
\begin{center}
\centerline{\psfig{figure=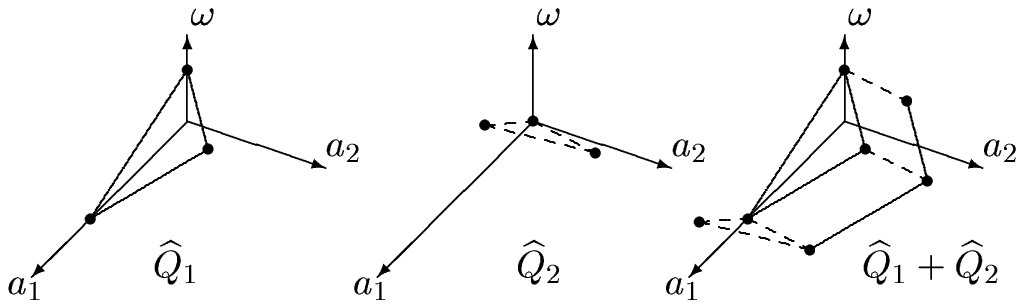}}
\caption{Lifted polytopes ${\widehat Q}_1$, ${\widehat Q}_2$, and
a regular mixed subdivision of ${\widehat Q}_1 + {\widehat Q}_2$.}
\label{fig3dmcc}
\end{center}
\end{figure}

To compute only the mixed volume,
the lift-and-prune approach was presented in~\cite{ec95},
using a primal model to prune in the tree of edge-edge combinations.
This approach operates in two stages.  First the polytopes are lifted
by adding one coordinate to every point in the supports.
In the second stage, one computes the facets of the lower hull of 
the Minkowski sum that are spanned by sums of edges.
These facets constitute the {\em mixed cells} in a mixed subdivision.
On the supports~$\cal A$ in~(\ref{eqexpolysup}),
we consider the lifted supports
\begin{equation} \label{eqlifted}
  {\widehat {\cal A}} = ({\widehat A}_1,{\widehat A}_2) \quad
   \begin{array}{l}
      {\widehat A}_1 = \{ (0,0,1) , (1,1,0) , (4,0,0) \} \\
      {\widehat A}_2 = \{ (0,0,0) , (1,2,0) , (3,1,1) \}
   \end{array}                 
\end{equation}
The lower hulls of the lifted polytopes are displayed in
Figure~\ref{fig3dmcc}.

\medskip

The two cells that contribute to the mixed volume are identified by
inner normals $\alpha$ and $\beta$ that satisfy
systems of linear equations and inequalities:

\smallskip

\begin{equation} \label{eqdual}
  \begin{array}{l}
    \! \! \! \! \alpha = (0,0,1) \\
    \left\{
       \begin{array}{rcl}
          4 \alpha_1 = \alpha_1 + \alpha_2 & < & 1 \\
          \alpha_1 + 2 \alpha_2 = 0 & < & 3 \alpha_1 + \alpha_2 + 1
       \end{array}
    \right.
  \end{array}
  \quad \quad
  \begin{array}{l}
    \! \! \! \! \beta = (2,-1,1) \\
    \left\{
       \begin{array}{rcl}
          \beta_1 + \beta_2 = 1 & < & 4 \beta_1 \\
          \beta_1 + 2 \beta_2 = 0 & < & 3 \beta_1 + \beta_2 + 1
       \end{array}
    \right.
  \end{array}
\end{equation}

\smallskip

\noindent These systems express that the cells correspond to facets 
spanned by the sum of two edges on the lower hulls of
${\widehat Q}_1$ and ${\widehat Q}_2$ respectively.
The lift-and-prune method with a dual version of the linear inequality
constrains as in~(\ref{eqdual}) was elaborated in~\cite{vgc96},
exploiting the fact that several polynomials can share the same
Newton polytope (see~\cite{hs95}) and with dimension reductions.

\medskip

\begin{figure}[hbt]
\begin{center}
\centerline{\psfig{figure=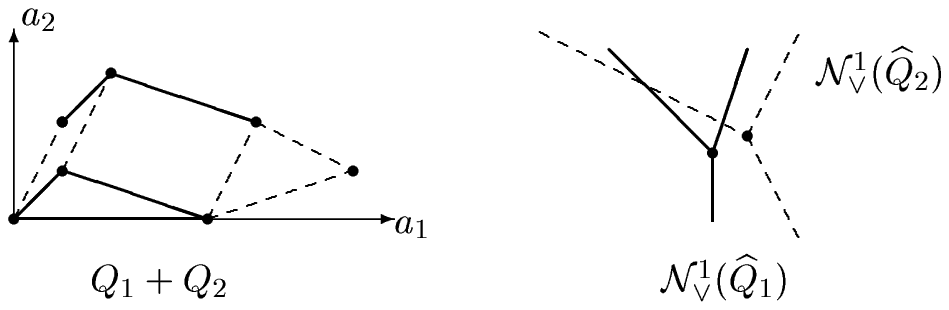}}
\caption{On the left we see the projection of a regular mixed subdivision 
of ${\widehat Q}_1 + {\widehat Q}_2$.  On the right, we have the dual 
construction with complexes ${\cal N}_{\vee}^{1}(Q_1)$ and 
${\cal N}_{\vee}^{1}(Q_2)$
collecting the cones of all vectors normal to the edges on the lower hulls
of ${\widehat Q}_1$ and ${\widehat Q}_2$ respectively.
The intersection of the cones contain the normals to the mixed cells.}
\label{figfans}
\end{center}
\end{figure}

The geometric dual construction to Figure~\ref{fig3dmcc} is
displayed in Figure~\ref{figfans}.

As in~\cite{hs95}, we assume that there are $r$ different Newton polytopes.
Given a tuple of lifted point sets
${\widehat {\cal A}}
 = ( {\widehat A}_1, {\widehat A}_2, \ldots, {\widehat A}_r )$,
any lifted cell ${\widehat {\cal C}}_{\bf v}$ of a regular subdivision 
can be characterized by its inner normal as
\begin{equation}
  {\widehat {\cal C}}_{\bf v} = (
 \partial_{\bf v} {\widehat A}_1, \partial_{\bf v} {\widehat A}_2, \ldots ,
 \partial_{\bf v} {\widehat A}_r ).
\end{equation}
Since ${\rm conv}({\widehat C}_{\bf v})
 = {\rm conv}(\sum_{i=1}^r \partial_{\bf v} {\widehat A}_i)$ 
is a facet of the lower hull,
the inner product $\langle . , {\bf v} \rangle$ attains \\[1mm]
its minimum over ${\widehat A}_i$ at $\partial_{\bf v} {\widehat A}_i$, i.e.,
\begin{equation} \label{equa}
 \forall {\widehat {\bf a}}, {\widehat {\bf b}}
        \in \partial_{\bf v} {\widehat A}_i: \
 \langle {\widehat {\bf a}} , {\bf v} \rangle
 = \langle {\widehat {\bf b}} , {\bf v} \rangle, \quad i = 1,2,\ldots,r,
\end{equation}
\begin{equation} \label{inequa}
 \forall {\widehat {\bf a}} \in {\widehat A}_i \setminus
         \partial_{\bf v} {\widehat A}_i, \
 \forall {\widehat {\bf b}} \in \partial_{\bf v} {\widehat A}_i: \
 \langle {\widehat {\bf a}} , {\bf v} \rangle
 > \langle {\widehat {\bf b}} , {\bf v} \rangle, \quad i = 1,2,\ldots,r.
\end{equation}
Algorithm~\ref{algshafac} (presented in~\cite{vgc96})
gives a way to compute all mixed cells by searching for feasible solutions
to the constraints~(\ref{equa}) and~(\ref{inequa}).
The algorithm generates a tree of all possible combinations of $k_i$-faces,
with feasibility tests to prune branches that do not lead to mixed cells.
The order of enumeration is organized so that mixed cells which
share some faces also share a part of the factorization work to be
done to solve the system defined by~(\ref{equa}).

\begin{algorithm}  \label{algshafac}
{\rm Pruning algorithm with shared factorizations subject to 
inequality constraints:

\begin{center}
\begin{tabular}{lcr}
Input: $({\widehat A}_1,{\widehat A}_2, \ldots, {\widehat A}_r)$,
  & & {\em lifted point sets} \ \\
 \ \ \ \ \ \ \ \ \ $\!$ ${\bf k} = (k_1,k_2,\ldots,k_r)$,
           $n = \sum_{i=1}^r k_i$, & & {\em $A_i$ appears $k_i$ times } \ \\
 \ \ \ \ \ \ \ \ \ $\!$
 $({\widehat {\cal F}}_1,{\widehat {\cal F}}_2, \ldots,
 {\widehat {\cal F}}_r)$.
  & \ \ & {\em $k_i$-faces of lower hull of ${\rm conv}({\widehat A}_i)$} \ \\ 
Output: ${\widehat {\frak S}_\omega} 
           = \{ \ {\widehat {\cal C}} \in {\widehat S}_\omega \ |
                \ V_n({\cal C},{\bf k}) > 0 \ \}$.
            & & {\em collection of lifted mixed cells} \ \\
\\
{\bf At level} $i$, $1 \leq i < r$: & & \\
\ \ {\em DATA and INVARIANT CONDITIONS}: & & \\
\ \ \ \ \ $\!$ $(M_1,\kappa)$: \ \
        $M_1 {\bf v} = {\bf 0} \not\Rightarrow v_{n+1} = 0$,
        $\kappa = {\displaystyle \sum_{j=1}^{i-1} k_j}$
                  & & \begin{tabular}{r}
                             {\em equalities} (\ref{equa}) \\
                             {\em upper triangular up to row $\kappa$}
                      \end{tabular} \\
\ \ \ \ \begin{tabular}{ll}
            $(M_2,\kappa)$: &
            $M_2 {\bf v} \geq {\bf 0} \not\Rightarrow -v_{n+1} \geq 0$ \\
            & $\dim(M_2) = n-\kappa$
          \end{tabular}
                 & & \begin{tabular}{r}
                             {\em inequalities} (\ref{inequa}) \\
                             {\em still feasible and reduced}
                      \end{tabular} \\
\ \ {\em ALGORITHM}: & & \\
\ \ \ \ \ for each ${\widehat C}_i \in {\widehat {\cal F}}_i$ loop
                         & & {\em enumerate over all $k_i$-faces} \ \\
\ \ \ \ \ \ \ {\em Triangulate}$(M_1,\kappa,{\widehat C}_i)$; 
                                   & & {\em ensure invariant conditions} \ \\
\ \ \ \ \ \ \ if $M_1 {\bf v} = {\bf 0} \not\Rightarrow v_{n+1} = 0$ 
                     & & {\em test for feasibility w.r.t.} (\ref{equa}) \ \\
\ \ \ \ \ \ \ \ then {\em Eliminate}$(M_1,M_2,\kappa,{\widehat C}_i,
                                                 {\widehat A}_i)$;
                & & {\em eliminate unknowns} \ \\
\ \ \ \ \ \ \ \ \ \ \ \ \ \ \ if $M_2 {\bf v} \geq {\bf 0} \not\Rightarrow 
                       -v_{n+1} \geq 0$ 
                & & {\em test for feasibility w.r.t.} (\ref{inequa}) \ \\
\ \ \ \ \ \ \ \ \ \ \ \ \ \ \ \ then proceed to next level $i+1$; & & \\
\ \ \ \ \ \ \ \ \ \ \ \ \ \ \ end if; & & \\
\ \ \ \ \ \ \ end if; & & \\
\ \ \ \ \ end for. & & \\
{\bf At level} $i=r$: & & \\
\ \ Compute $\bf v$: $M_1 {\bf v} = {\bf 0}$; & & \\
\ \ Merge the new cell ${\cal C}_{\bf v}$ with the 
      list ${\widehat {\frak S}_\omega}$. & &
\end{tabular}
\end{center}
}
\end{algorithm}
Note that (\ref{inequa}) has to be weakened to $\geq$ type inequalities
in order to be able to compute also subdivisions that are not fine.
This also explains the merge operation at the end.
The feasibility tests in the algorithm allow an efficient computation
of the mixed cells.
The conditions~(\ref{equa}) and~(\ref{inequa}) are verified incrementally.
After choosing a $k_i$-face ${\widehat C}_i = \{
{\widehat {\bf c}}_{0i}, {\widehat {\bf c}}_{1i} , \ldots ,
{\widehat {\bf c}}_{k_i i} \}$ of ${\widehat A}_i$, linear programming
is used to check whether $({\widehat C}_1, \ldots , {\widehat C}_i )$
can lead to a mixed cell in the induced subdivision.

\medskip

We end this section with complexity results.
The complexity of computing mixed volumes is proven~\cite{dgh98}
to be $\#P$-hard.  This complexity class is typical for all enumerative
problems, since, unlike the class $NP$, there exists no algorithm that
runs in polynomial time for arbitrary dimensions to verify that
a guessed answer is correct.
Although the current algorithmical practice suggests that computing mixed 
volumes is harder than computing volumes of polytopes (which is also
known as a $\#P$-hard problem~\cite{df88}), this is not the case
from a complexity point of view as shown in~\cite{dgh98}.
In~\cite{emi96} it is shown that the mixed volume~$V_n({\cal Q})$ 
is bounded from below by $n! {\rm vol}_n(Q_\mu)$, $Q_\mu$ being the
polytope of minimum volume in~$\cal Q$.

\subsection{Sparse Polynomial Systems solved by Polyhedral Homotopies}
\label{secsparhom}

The simplest system in the polytope model that still has isolated solutions
in~$(\cc^*)^n$ has exactly two terms in every equation.
Polyhedral homotopies~\cite{hs95} solve systems with random complex
coefficients starting from sparser subsystems.  For~(\ref{eqexsys}),
the homotopy with supports~$\widehat {\cal A}$ as in~(\ref{eqlifted}) is

\begin{equation} \label{eqexpolyhom}
  {\widehat P}(x_1,x_2,t) =
    \left\{
      \begin{array}{r}
         c_1 x_1^4 t^0 + c_2 x_1 x_2 t^0 + c_3 t^1 = 0 \\
         c'_1 x_1^3 x_2 t^1 + c'_2 x_1 x_2^2 t^0 + c'_3 t^0 = 0
      \end{array}
    \right.
  \quad {\rm with~} c_i,c'_i \in \cc^*.
\end{equation}
The exponents of~$t$ are the values of the lifting~$\omega$
applied to the supports.

\medskip 

To find the start systems, we look at Figure~\ref{fig3dmcc},
at the subdivision that is induced by this lifting process.
The start systems have Newton polytopes spanned by one edge of the first
and one edge of the second polytope.
Since the two cells that contribute to the mixed volume
are characterized by their inner normals $\alpha$ and $\beta$
satisfying~(\ref{eqdual}) we denote the start systems respectively
by $P^\alpha$ and $P^\beta$.
To compute start solutions, unimodular transformations make the
system triangular as follows.
After dividing the equations so that the constant term is present,
we apply the substitution $x_1 = y_2$, $x_2 = y_1^{-1} y_2^3$ on
$P^\alpha$ as follows:

\begin{equation} \label{eqbin}
  P^{\alpha}({\bf x})
  = \left\{
      \begin{array}{r}
         x_1^3 x_2^{-1} + c''_1 = 0 \\
         x_1 x_2^2 + c''_2 = 0
      \end{array}
    \right.
\quad
  P^{\alpha}({\bf x} = {\bf y}^U)
  = \left\{
      \begin{array}{r}
         y_1 + c''_1 = 0 \\
         y_1^{-2} y_2^7 + c''_2 = 0
      \end{array}
    \right.
\end{equation}

\noindent The substitution in (\ref{eqbin}) is apparent in the notation
(as used in~\cite{li97})
${\bf x}^V = ({\bf y}^U)^V = {\bf y}^{VU} = {\bf y}^L$ elaborated as
\begin{equation}
  \begin{array}{rcl}
    \left(
      \begin{array}{c}
         x_1^3 \cdot x_2^{-1} \\ x_1^1 \cdot x_2^2
      \end{array}
    \right) 
    & = &
    \left(
      \begin{array}{c}
         (y_1^0 y_2^1)^3 \cdot (y_1^{-1} y_2^3)^{-1} \\
         (y_1^0 y_2^1)^1 \cdot (y_1^{-1} y_2^3)^2
      \end{array}
    \right) \\
    & = &
    \left(
      \begin{array}{c}
         y_1^{3 \cdot 0 - 1 \cdot (-1)} \cdot y_2^{3 \cdot 1 - 1 \cdot 3} \\
         y_1^{1 \cdot 0 + 2 \cdot (-1)} \cdot y_2^{1 \cdot 1 + 2 \cdot 3}
      \end{array}
    \right)
    = 
    \left(
      \begin{array}{c}
         y_1^1 \cdot y_2^0 \\ y_1^{-2} \cdot y_2^7
      \end{array}
    \right).
  \end{array}
\end{equation}
The exponents are calculated by the factorization $VU = L$:
\begin{equation}
  \left[ \begin{array}{rr}
            3 & -1 \\ 1 & 2
         \end{array}
  \right]
  \left[ \begin{array}{rr}
            0 & 1 \\ -1 & 3
         \end{array}
  \right] =
  \left[ \begin{array}{rr}
            1 & 0 \\ -2 & 7
         \end{array}
  \right] .
\end{equation}
Since $\det(U) = 1$, the matrix $U$ is called unimodular.

\medskip

The polyhedral homotopy~(\ref{eqexpolyhom}) directly extends
the solutions of $P^{\alpha}$ to the target system.
To obtain a homotopy starting
at the solutions of~$P^{\beta}$, we substitute in~(\ref{eqexpolyhom})
$x_1 \leftarrow {\widetilde x}_1 t^{\beta_1}$, 
$x_2 \leftarrow {\widetilde x}_2 t^{\beta_2}$
and clear out the lowest powers of~$t$.
This construction appeared in~\cite{hs95} and provides an algorithmic
proof of the following theorem.

\begin{theorem} \label{theobera}
(Bernshte\v{\i}n~{\rm \cite[Theorem A]{ber75}})
For a general choice of coefficients for~$P$,
the system $P({\bf x}) = {\bf 0}$ has
exactly as many regular solutions as its mixed volume~$V_n({\cal Q})$.
\end{theorem}
The original algorithm Bernshte\v{\i}n used in his proof was implemented
in~\cite{vvc94}.

\medskip

For the numerical stability of polyhedral continuation, it is important
to have subdivisions induced by low lifting values, since those influence
the power of the continuation parameter.
In~\cite{vgc96} explicit lower bounds on integer lifting values were derived,
but unfortunately the dynamic lifting algorithm does not generalize that
well~\cite{lw98} if one is only interested in the mixed cells of a mixed
subdivision.  A balancing method was proposed in~\cite{glvw98} to improve
the stability of homotopies induced by random floating-point lifting values.

\medskip

Once all solutions to a polynomial system with randomly generated coefficients
are computed, we use cheater's homotopy to solve any specific system with
the same Newton polytopes. 
One could say that polyhedral homotopies have removed the cheating part.
The main advantage of polyhedral methods is that the mixed volume is a much
sharper root count in most applications, leading to fewer paths to trace.
They also allow more flexibility to exploit symmetry as demonstrated
in~\cite{vg95}.

\medskip

In case the system has fewer isolated solutions than the mixed volume,
we consider the face systems.  Define the face of a
polynomial $p$ with support $A$ as follows:
\begin{equation}
  p({\bf x}) = \sum_{{\bf a} \in A} c_{\bf a} {\bf x}^{\bf a}
  \quad \mbox{has faces} \quad
  \partial_{\bf v} p({\bf x})
  = \sum_{{\bf a} \in \partial_{\bf v} A} c_{\bf a} {\bf x}^{\bf a}
  \quad \mbox{for } {\bf v} \not= {\bf 0}.
\end{equation}
For ${\bf v} \not= {\bf 0}$, the corresponding face system of $P$ 
is $\partial_{\bf v} P = ( \partial_{\bf v} p_1,
\partial_{\bf v} p_2, \ldots , \partial_{\bf v} p_n )$.

\begin{theorem} \label{theoberb}
(Bernshte\v{\i}n~{\rm \cite[Theorem B]{ber75}})
Suppose $V_n({\cal Q}) > 0$.  Then,
$P({\bf x}) = {\bf 0}$ has fewer than~$V_n({\cal Q})$
isolated solutions if and only if $\partial_{\bf v} P({\bf x}) = {\bf 0}$
has a solution in $(\cc^*)^n$, for~${\bf v} \not= {\bf 0}$.
\end{theorem}
As is the case for our running example~(\ref{eqexsys}),
the Newton polytopes may be in generic position such that for any
nonzero choice of the coefficients, the system has exactly as
many isolated solutions as the mixed volume.
In practical applications however, 
how can we decide whether paths are really going towards infinity?
Relying on the actual computed values is arbitrarily,
because $10^4$ is as far from infinity as~$10^8$, so we need
algebraic structural data to certify the divergence.

\medskip

In the polyhedral end game~\cite{hv98} solution paths are represented
by power series expansions:
\begin{equation} \label{eqpower}
  \left\{
    \begin{array}{rcl}
       x_i(s) & = & a_i s^{v_i} ( 1 + O(s) ) \\
         t(s) & = & 1 - s^m
    \end{array}
  \right.
  \quad \quad t \approx 1, \quad s \approx 0.
\end{equation}
The winding number $m$ is lower than or equal to the multiplicity of the
solution.  For a solution diverging to infinity or to a zero-component
solution we observe that $v_i \not= 0$.
According to Theorem~\ref{theoberb}, this solution vanishes at a
face system $\partial_{\bf v} P$
(same ${\bf v}$ with components $v_i$ as in~(\ref{eqpower})),
certifying the divergence.

\medskip

To check whether a solution path really diverges is equivalent to
the test on the value for~$v_i$.  A first-order approximation
of $v_i$ can be computed by 
\begin{equation}
   \frac{\log|x_i(s_1)| - \log|x_i(s_0)|}{\log(s_1) - \log(s_0)}
   = v_i + O(s_0),
\end{equation}
with $0 < s_1 < s_0$.  The above formula assumes the correct value for~$m$.
To compute $m$, solution paths are sampled geometrically with ratio $h$ 
as $s_k = h^{k/m} s_0$.  The errors on the estimates for~$v_i$ are
\begin{eqnarray}
  e_i^{(k)} & = & (\log|x_i(s_k)|     - \log|x_i(s_{k+1})|)
                - (\log|x_i(s_{k+1})| - \log|x_i(s_{k+2})|) \\
            & = & c_1 h^{k/m} s_0 (1 + O(h^{k/m})).
\end{eqnarray}
An estimate for~$m$ is derived from two consecutive errors~$e_i^{(k)}$.
Extrapolation improves this estimate.
So, by an inexpensive side calculation at the end of the paths, we obtain
important structural algebraic information about the system.

\medskip

Recall that $V_n({\cal Q})$ count the roots in~$(\cc^*)^n$.
Using Newton polytopes to count affine roots
(i.e.: in~$\cc^n$ instead of~$(\cc^*)^n$)
was proposed in~\cite{roj94} with the notion of shadowed polytopes obtained by
the substitution $x_i \leftarrow x_i + c_i$ for arbitrary constants $c_i$.
To arrive at sharper bounds, it suffices (see~\cite{lw96} and~\cite{rw96})
to add a random constant to every equation.
Stable mixed volumes~\cite{hs97} provide a generically sharp affine root count.
The constructions in~\cite{ev99} and~\cite{glw98} avoid the use of 
recursive liftings to compute stable mixed volumes.
Further developments and generalizations can be found in~\cite{roj99}.

\subsection{Determinantal Polynomials arising in Enumerative Geometry}

Homotopies for solving problems in enumerative geometry appeared
in~\cite{hss98}.  The algorithms in the numerical Schubert calculus
originated from questions in real enumerative geometry~\cite{sot97a,sot97b}
and have their main application to the pole placement problem~\cite{byr89},
\cite{rrw96,rrw98}, \cite{ros94}, \cite{rw99} in control theory.

\medskip

The enumerative problems are formalized in some ``finiteness'' theorems,
avoiding the explicit but involved (as in~(\ref{eqgrassdeg})) formulas
for the root counts.

\begin{theorem} \label{theoprob1}
The number of $r$-planes meeting $mr$ general $m$-planes in~$\cc^{m+r}$
is a finite constant.
\end{theorem} 

The first homotopy presented in~\cite{hss98} uses a Gr\"obner basis
for the ideal that defines $G_{mr}$, as is derived in~\cite{stu93}.
By Gr\"obner bases questions concerning any polynomial system are solved
by relation to monomial equations.
Every Gr\"obner basis defines a flat deformation, which preserves
the structure of the solution set~\cite{eis95}.
Geometrically, this type of deformation is used to collapse the solution
set in projective space to the coordinate hyperplanes, or in the opposite
direction, to extend the solutions of the monomial equations to those
of the original system.                                       
The flat deformations that are obtained in this way are similar to 
toric deformations in the sense that one moves from the solutions of
a subsystem to the solutions of the whole system.

\medskip

The Gr\"obner homotopies of~\cite{hss98} work in the synthetic 
Pl\"ucker embedding, and need to take the large set of defining 
equations of $G_{mr}$ into account.
When expanding the minors into local coordinates, these equations
are automatically satisfied, which leads to a much smaller
polynomial system.  Consequently, the second type of homotopies
of~\cite{hss98}, the so-called SAGBI homotopies are more efficient.
Instead of an ideal, we now have a subalgebra and work with
a SAGBI basis,
i.e.: the Subalgebra Analogue to a Gr\"obner Basis for an Ideal. 
The term order selects the monomials on the diagonal as the dominant
ones.  This implies that in the flat deformation (see~\cite{stu96}
for a general description) only the diagonal monomials remain at~$t=0$.

\medskip

For $m=2=r$, the equations of the SAGBI homotopy in determinantal form are
\begin{equation} \label{eqsabi22det}
  p_i({\bf x}) =
  \det \left[
    \left.
    \begin{array}{cc}
       c_{11}^{(i)} & c_{12}^{(i)} \\
       c_{21}^{(i)} & c_{22}^{(i)} \\
       c_{31}^{(i)} & c_{32}^{(i)} \\
       c_{41}^{(i)} & c_{42}^{(i)} \\
    \end{array}
    \right|
    \begin{array}{cc}
       x_{11} ~ & x_{12} \\
       x_{21} t & x_{22} \\
              1 & 0 \\
              0 & 1 \\
     \end{array}
  \right]
  = 0, \quad i=1,2,3,4,
\end{equation}
where the coefficients $c_{kl}^{(i)}$ are random complex constants. 
In expanding the minors of~(\ref{eqsabi22det}), the lowest power of $t$
is divided out, minor per minor.
The system at $t=0$ is solved by polyhedral continuation.
The system at $t=1$ serves as start system in the cheater's homotopy
to solve any system with particular real values for the
coefficients~$c_{kl}^{(i)}$.
Figure~\ref{figsagbi} outlines the structure of the general solver.

{\small
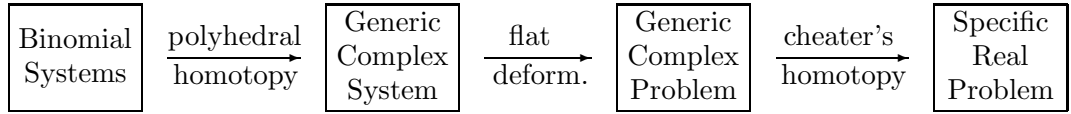
\begin{figure}[hbt]
\begin{center}
\begin{picture}(400,40)(0,0)

\put(0,0){\framebox(50,40)[c]
              {\begin{tabular}{c} Binomial \\Systems \end{tabular}}} 

\put(60,20){\vector(1,0){50}}  \put(60,25){polyhedral}
                               \put(62,10){homotopy}

\put(120,0){\framebox(50,40)[c]
              {\begin{tabular}{c} Generic \\ Complex \\ System \end{tabular}}} 

\put(180,20){\vector(1,0){40}}  \put(189,25){flat}
                                \put(183,10){deform.}

\put(230,0){\framebox(50,40)[c]
              {\begin{tabular}{c} Generic \\ Complex \\ Problem \end{tabular}}} 

\put(290,20){\vector(1,0){50}}  \put(293,25){cheater's}
                                \put(292,10){homotopy}

\put(350,0){\framebox(50,40)[c]
              {\begin{tabular}{c} Specific \\ Real \\ Problem \end{tabular}}} 

\end{picture}
\caption{The SAGBI homotopy is at the center of the concatenation.}
\label{figsagbi}
\end{center}
\end{figure}
}

SAGBI homotopies have been implemented~\cite{ver98} to
verify some large instances of input planes for which it was conjectured
that all solution planes would be real.
We refer to~\cite{rs98} and \cite{sot98} for related work on these conjectures.
In~\cite{sot99b} an asymptotic choice of inputs is generated for which all
solutions are proven to be real.

\medskip

The third type of homotopies presented in~\cite{hss98} are the
so-called Pieri homotopies.
Since they are closer to the intrinsic geometric nature, they are
applicable to a broader class of problems.  In particular, we obtain
an effective proof for the following.

\begin{theorem} \label{theoprob2}
The number of $r$-planes meeting $n$ general $(m+1-k_i)$-planes in~$\cc^{m+r}$,
with $k_1 + k_2 + \cdots + k_n = mr$, is a finite constant.
\end{theorem}
Note that when all $k_i = 1$, we arrive at Theorem~\ref{theoprob1}.
For general $k_i \not= 1$, we are not aware of any explicit formulas
for the number of roots.

\medskip

Figure~\ref{figcelldeco} shows a part of a cellular decomposition
of $G_{22}$ with the determinantal equations.
We specialize the pattern~$X$ that represents a solution line by
setting some of its coordinates to zero.
This specialization determines a specialization of the input lines
as follows: take those basis vectors not indexed by rows of~$X$
where zeroes have been introduced.  The special line $S_X$ for this
example is as in~(\ref{eqpierihom}) spanned by the first and third
standard basis vector.

\begin{figure}[hbt]
\begin{center}
\begin{picture}(240,120)(20,-10)

\thicklines

\put(-50,85){\small
${\rm det}\left[ 
  S_X
  \left|
     \begin{array}{cc}
        x_{11} &   0    \\
          0    &   0    \\
          0    & x_{32} \\
          0    & x_{42} \\
     \end{array}
  \right.
 \right] = 0$}

\put(140,85){\small
${\rm det}\left[ 
  S_X
  \left|
     \begin{array}{cc}
        x_{11} &   0    \\
        x_{21} &   0    \\
          0    & x_{32} \\
          0    &   0    \\
     \end{array}
  \right.
 \right] = 0$}

\put(110,50){\line(4,1){40}}
\put(90,50){\line(-4,1){40}}

\put(10,15){\small
${\rm det}[L_3 | X]
 = {\rm det}\left[ 
   L_3
   \left|
   \begin{array}{cc}
      x_{11} &   0    \\
      x_{21} &   0    \\
        0    & x_{32} \\
        0    & x_{42} \\
   \end{array}
  \right.
 \right] = 0$}

\put(260,10){[2~4]}
\put(227,30){[1~4]}
\put(293,30){[2~3]}

\put(257,18){\line(-2,1){10}}
\put(285,18){\line(2,1){10}}

\end{picture}
\caption{Part of a cellular decomposition of the Grassmannian of all 2-planes.
At the right we have the short-hand notation with brackets.
The bracket~$[2~4]$ contains the row indices to the lowest nonzero entries
in~$X$.  }
\label{figcelldeco}
\end{center}
\end{figure}

\medskip

Figure~\ref{figcelldeco} pictures patterns of the moving 2-planes in the Pieri
homotopy algorithm for the case $(m,r) = (2,2)$, see Figure~\ref{figlines}.
The bottom matrix is the general representation of a solution that
intersects already the two input lines spanned by the standard basis
vectors.  At the leaves of the tree by linear algebra operations we can 
intersect with a third input line.  Moving down the poset, we deform form 
the left configuration in Figure~\ref{figlines} to the general problem.

\medskip

Denote by $L_1$ and $L_2$ the lines already met by $X$.
At the leaves of the tree in Figure~\ref{figcelldeco}
we intersect with the fourth line $L_4$.
The special position for the third line $L_3$ is represented by the 
matrix $S_X$, which intersects any $X$ with coordinates as at the leaves
of the tree.
In the homotopy $H(X,t) = {\bf 0}$ we deform the line spanned by the columns
of $S_X$ to line $L_3$, for $t = 0$ to $t = 1$.
\begin{equation} \label{eqpierihom}
   S_X = \left[
          \begin{array}{cc}
             1 & 0 \\
             0 & 0 \\
             0 & 1 \\
             0 & 0 \\
          \end{array}
       \right]
   \quad \quad
   H(X,t) =
     \left\{
        \begin{array}{l}
           \det( L_4 | X ) = 0 \\
           \det( (1-t) S_X + tL_3 | X ) = 0 \\
        \end{array}
     \right.
   \quad t \in [0,1].
\end{equation}
Every solution $X(t)$ of $H(X(t),t) = {\bf 0}$ intersects already three
lines: $L_1$, $L_2$ and $L_3$.
At the end, for $t=1$, $X$ also meets the line~$L_3$ in a point.

\medskip

The homotopy~(\ref{eqpierihom}) deforms two solution lines,
starting at patterns which have their row indices for the lowest nonzero
entries respectively as in~$[1~4]$ and in~$[2~3]$.
The correctness of this homotopy (proven in~\cite{hss98} and~\cite{hv99})
justifies the formal root count using the localization poset.
This combinatorial root count proceeds in two stages.  First we build
up the poset, from the bottom up, diminishing the entries in the brackets
under the restriction that the same entry never occurs twice or more.
Secondly, we descend from the top of the poset, collecting and adding
up the counts at the nodes in the poset.
More examples and variations are in~\cite{hv99}.

\medskip

To solve the general intersection problem of Theorem~\ref{theoprob2}, the
special $(m+1-k_i)$-planes lie in the intersection of special $m$-planes.
In the construction of the poset one has to follow additional rules
as to ensure a solution that meets the intersection of special $m$-planes.
We refer to~\cite{hv99} for details.

\medskip

The third enumerative problem we can solve is formalized as follows.

\begin{theorem} \label{theoprob3}
The number of all maps of degree $q$ that produce $r$-planes in $\cc^{m+r}$
meeting $mr + q(m+r)$ general $m$-planes at specified interpolation
points is a finite constant.
\end{theorem}
In~\cite{rrw96,rrw98} and~\cite{wrr94} explicit formulas are given 
for this finite constant along with other combinatorial identities.
Following a hint of Frank Sottile (see also~\cite{sot99a})
and reverse engineering on the root counts in~\cite{rrw96},
Pieri homotopies were developed in~\cite{hv99} whose correctness
yields a proof for Theorem~\ref{theoprob3}.

\medskip

The analogue to Figure~\ref{figcelldeco} for maps of degree one
into $G_{22}$ is displayed in Figure~\ref{figcurvcelldeco}.

\begin{figure}[hbt]
\begin{center}
\begin{picture}(220,120)(20,-10)

\thicklines

\put(-65,85){\small
${\rm det}\left[ 
  S_X
  \left|
    \begin{array}{cc}
      x^0_{11} & x^1_{12} s~~~ \\
      x^0_{21} & ~~~x^0_{22} t \\
         0     & ~~~x^0_{32} t \\
         0     & ~~~x^0_{42} t \\
    \end{array}
  \right.
 \right] = 0$}

\put(140,85){\small
${\rm det}\left[ 
  S_X
  \left|
   \begin{array}{cc}
      x^0_{11} & 0 \\
      x^0_{21} & x^0_{22}  \\
      x^0_{31} & x^0_{32} \\
         0     & x^0_{42} \\
   \end{array}
  \right.
 \right] = 0$}

\put(110,50){\line(4,1){40}}
\put(90,50){\line(-4,1){40}}

\put(-10,15){\small
${\rm det}[L_n | X(s,t)]
 = {\rm det}\left[ 
   L_n
   \left|
   \begin{array}{cc}
      x^0_{11} & x^1_{12} s~~~ \\
      x^0_{21} & ~~~x^0_{22} t \\
      x^0_{31} & ~~~x^0_{32} t \\
         0     & ~~~x^0_{42} t \\
   \end{array}
  \right.
 \right] = 0$}

\put(260,10){[3~5]}
\put(227,30){[2~5]}
\put(293,30){[3~4]}

\put(257,18){\line(-2,1){10}}
\put(285,18){\line(2,1){10}} 

\end{picture}
\caption{Part of a cellular decomposition of the Grassmannian of 
maps of degree~1 that produce 2-planes in projective 3-space.
The bracket notation at the right corresponds to a matrix representation
of the coefficients of the map~$X(s,t)$.  }
\label{figcurvcelldeco}
\end{center}
\end{figure}

To solve the problem in Theorem~\ref{theoprob3} we need a special position
for the interpolation points.  By moving those to infinity, the
dominant monomials in the maps allow to re-use the same special $m$-planes,
whose entries should be considered modulo~$m+r$.
The homotopy to satisfy the $n$-th intersection condition is:

\begin{equation} \label{eqpieriqhom}
  H(X(s,t),s,t) =
  \left\{
    \begin{array}{rl}
       \det(L_i | X(s_i,t_i)) = 0 & i=1,2,\ldots,n-1 \\
       \det((1-t) S_X + t L_n | X(s,t) ) = 0 & \\
       (s-1)(1-t) + (s-s_n)t = 0 & t \in [0,1]
    \end{array}
   \right.
\end{equation}  
Note that the continuation parameter~$t$ moves the interpolation point 
from infinity, at $(s,t) = (1,0)$, to the specific value $(s,t) = (s_n,1)$.

\medskip

See~\cite{sot99a} for information on the selection of the input planes so
that all maps are real.

\medskip

As an example of another problem in enumerative geometry we mention
the 27 lines on a cubic surface in 3-space.
According to~\cite{mum76}, this is one of the gems
hidden in the rag-bag of projective geometry.
In~\cite{seg42}, the 27 lines are determined by breaking up the cubic
surface into three planes in a continuous way such that each intermediate
position is nonsingular.  It is shown that this continuous variation
is also valid in the real field.

\section{Numerical Software for Solving Polynomial Systems}

In computer algebra one wants to compute exactly as long as possible
and to defer the approximate calculations to the very last end.
Exactly the opposite way is taken in homotopy methods:
here we use floating-point arithmetic and only increase the precision
when needed.

\medskip

Next we mention programs with special features for polynomial systems.
See~\cite[Chapter VIII]{ag97} for a list of available software for path 
following.
HOMPACK~\cite{msw89,wbm87} is a general continuation package with a 
polynomial driver.  It has been parallelized~\cite{acw89,hw89},
extended with an end game~\cite{sws96}, and upgraded~\cite{wsmmw97}
to Fortran~90.  POLSYS\_PLP~\cite{wsw98} provides linear-product root
counts to be used in conjunction with HOMPACK90.
The Fortran code for CONSOL is contained in~\cite[Appendix 6]{mor87}.
The C-program pss~\cite{mal96} applies homotopy continuation with
verification by $\alpha$-theory.  Pelican~\cite{hub95,hub96}
implements in~C the polyhedral methods of~\cite{hs95}.
Efficient Fortran software for polyhedral continuation with facilities 
to compute all affine roots is used in~\cite{glw98}.
The computation of mixed volumes with the C program mvlp~\cite{emi94,ec95}
is a crucial step for sparse resultants~\cite{wem98}.
A distributed version has been created in~\cite{gio96}.

\medskip

PHC is written in Ada and originated during the doctoral research of 
the author~\cite{ver96}.  Executable versions were first released
at the PoSSo open workshop on software~\cite{ver95}.
The public release of the sources is described in~\cite{ver97}.
The package is organized as a tool box, organized along four stages of 
the solver.  Figure~\ref{figstages} presents the flow of the solver.
The package is menu-driven and file-oriented.
A general-purpose black-box solver is available. 

{\small
\begin{figure}[hbt]
\begin{center}
\begin{picture}(360,145)(10,0)

\put(20,120){\framebox(130,20)[c]{1. Preconditioning}}
\put(50,100){$\diamond$ Coefficient Scaling}
\put(50,85){$\diamond$ Reduction of degrees}

\put(205,120){\framebox(130,20)[c]{4. Validation}}
\put(235,100){$\diamond$ Refining of the roots}
\put(235,85){$\diamond$ Analysis of condition $\#$s}

\put(20,50){\framebox(130,20)[c]{2. Root Counting}}
\put(30,35){$\diamond$ B\'ezout : degrees}
\put(30,20){$\diamond$ Bernshtein : polytopes}
\put(30,5){$\diamond$ Schubert : SAGBI/Pieri}

\put(205,50){\framebox(140,20)[c]{3. Homotopy Continuation}}
\put(215,30){$\diamond$ Fix continuation parameters}
\put(215,15){$\diamond$ Choose Predictor-Corrector}

\put(160,57){\line(1,0){25}}
\put(160,60){\line(1,0){25}}
\put(185,64){\line(2,-1){11}}
\put(185,53){\line(2,1){11}}

\put(34,113){\line(0,-1){25}}
\put(38,113){\line(0,-1){25}}
\put(31,91){\line(1,-2){5}}
\put(41,91){\line(-1,-2){5}}

\put(219,107){\line(0,-1){25}}
\put(223,107){\line(0,-1){25}}
\put(216,102){\line(1,2){5}}
\put(226,102){\line(-1,2){5}}

\end{picture}
\caption{The four stages in the flow of the PHC solver.}
\label{figstages}
\end{center}
\end{figure}
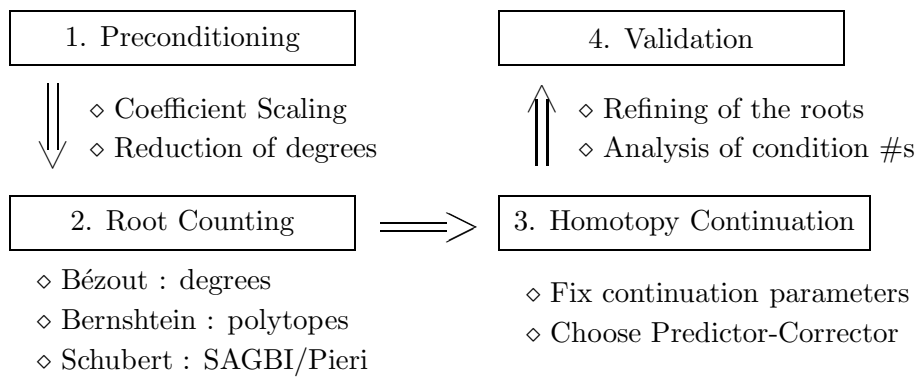
}

The new second release of PHC uses Ada~95 concepts in the construction
of the mathematical library.
It is developed with the freely available gnu-ada compiler
(currently at version 3.11p) on various platforms.
To run the software no compilation is needed,
as binaries are available for Unix Workstations: running SUN Solaris and
SGI Irix, and Pentium PCs: running Linux and Solaris.
The portability of PHC is ensured by the gnu-ada compiler.

\medskip

Another main feature of the second release are the homotopy methods
for the Schubert calculus.  Implementing those homotopies was a matter
of plugging in the equations and calling the path trackers.
The third release of the package should offer a more comprehensive
environment to construct homotopies, providing an easier access to
the two main computational engines: mixed-volume computation and
polynomial continuation.

\section{The Database of Applications}

The polynomial systems in scientific and engineering models
are a continuing source of open problems.
Systems that come from academic questions are often conjectures providing
computational evidence in a developing theory.
In various engineering disciplines polynomial systems represent a
modeling problem, e.g.: a mechanical device.
The origin of a polynomial system matters when the original problem
formulation does not admit well-conditioned solutions.
As a general method to deal with badly scaled systems to compute equilibria
of chemical reaction systems, coefficient and equation scaling
was developed in~\cite{mm87},
see also~\cite[Chapter 5]{mor87} and~\cite{wbm87}.

\medskip

The collection of test systems is organized as a database and
available via the author's web pages,
A good test example reveals properties of the solution method and
has a meaningful application.
Besides the algebraic formulation it contains the fields:
title (meaningful description), references (problem source), root counts
(B\'ezout bounds and mixed volume), and solution list.

\medskip

\noindent Instead of producing a huge list with an overview, we pick
some important case studies.

\begin{description}

\item[{\bf katsura-n}]  ({\em magnetism problem~{\rm \cite{kat94}}})
  The number of solutions equals the total degree $D = 2^n$,
  so the homotopy based on~$D$ is optimal to solve this problem.
  Because the constant term is missing in all except one equation, the system
  is an interesting test problem for affine polyhedral methods.

\smallskip

\item[{\bf camera1s}] ({\em computer vision~{\rm \cite{fm90}}})
  The system models the displacement of a camera between two positions
  in a static environment~\cite{emi94}.
  The multi-homogeneous homotopy is optimal for this problem, requiring
  20 solution paths to trace instead of~$D = 64$.

\smallskip

\item[{\bf game$n$two}] ({\em totally mixed Nash equilibria for $n$ players
  with two strategies~{\rm \cite{mm97,mcl97}}})
  This is another instance where multi-homogeneous homotopies are optimal.
  The number of solutions grows like $n! e^{-1}$ as $n \rightarrow \infty$.
  The largest system that is currently completely solvable is for $n = 8$
  requiring 14,833 paths to trace.
  Situations exist for which all solutions are meaningful.

\smallskip

\item[{\bf cassou}] ({\em real algebraic geometry})
  This system illustrates the success story of polyhedral homotopies:
  the total degree equals~1,344, best known B\'ezout bound is 312
  (see~\cite{lww96}), whereas the mixed volume gives~24.
  Still eight paths are diverging to infinity and polyhedral end
  games~\cite{hv98} are needed to separate those diverging paths from the
  from the other finite ill-conditioned roots.

\smallskip

\item[{\bf cyclic-n}] ({\em Fourier transforms~{\rm \cite{bjo89,bf91}}})
  For $n=7$, polyhedral homotopies are optimal, with all 924 paths leading
  to finite solutions.  For $n \geq 8$, the mixed volume overestimates the
  number of roots and there are components of solutions.
  In~\cite{sv99} the degrees of the components were computed for $n=8,9$.
  There are 34940 cyclic 10-roots, generated by 1747 solutions.

\smallskip

\item[{\bf pole28sys}] ({\em pole placement problem~{\rm \cite{byr89}}})
  This system illustrates the efficiency of SAGBI homotopies for verifying
  a conjecture in real algebraic geometry~\cite{sot98}.
  With the input planes chosen to osculate a rational normal curve,
  an instance with all 1,430 solutions real and isolated was solved
  in~\cite{ver98}.  The problem is relevant to control theory~\cite{rw99}.

\smallskip

\item[{\bf stewgou40}] ({\em mechanism design~{\rm \cite{die98}}})
  Whether the Stewart-Gough parallel platform in robotics could have all
  its 40 solutions real was a notorious open problem, until recently,
  as it was solved by numerical continuation methods~\cite{die98}.
  The problem formulation in~\cite{die98} is highly deficient: the
  mixed volume equals 1,536 whereas only 40 solution paths will converge.

\end{description}
We emphasize that we have optimal homotopies for three classes
of polynomial systems, but not for all possible structures.
Although one can solve a modelling problem by a black-box polynomial-system
solver, knowing the origin of the problem leads in most cases to more favorable
algebraic formulations that help the resolution of a polynomial system.
To produce really meaningful solutions one often has to be close 
to the source of the problems and be able to interact with the people who
formulate the polynomial systems. 

\medskip

In closing this section we list some notable usages of PHC.
Charles Wampler~\cite{wam96} used a preliminary version of PHC to count
the roots of various systems in mechanical design.
Root counts for linear subspace intersections in the Schubert calculus
were computed by Frank Sottile, see~\cite{sot98} for various tables.
A third example comes from computer graphics.
To show that the 12 lines tangent to four given spheres can all be real,
Thorsten Theobald used PHC, choosing appropriate parameters in the algebraic
formulation set up by Cassiano Durand.

\section{Closing Remarks and Open Problems}

The three classes presented in this paper are by no means exhaustive, but
give an idea of what can be done with homotopies to solve polynomial systems.
The root counts constitute the theoretical backbone for general-purpose
black-box solving.  Yet, the homotopy methods are flexible enough to
exploit a particular geometrical situation, with guaranteed optimal
complexity when applied to generic instances.

\medskip

From algebraic geometry formal procedures based on intersection theory 
count the number of solutions to classes of polynomial systems.  
Examples are the theorems of B\'ezout, Bernshte\v{\i}n and Schubert.
For these situations we construct a start system and have a homotopy to deform
the solutions to this start system to the solutions to any specific problem.
There are many other cases for which one knows how to count but not how
to deform and solve efficiently.
Research in homotopy methods is aimed at turning the formal root counts
into effective numerical methods. 
As open problem we can ask for a meta-homotopy method to connect formal root 
counting methods to solving generic systems and deformation procedures.

\medskip

In most applications, only the real solutions are important.
Once we know an optimal homotopy to solve the problem in the complex case,
we would like to know whether all solutions can be real and how the real
solutions are distributed.  The reality question appears for instance in
the theory of totally mixed Nash equilibria and in the pole placement
problem.  Finding well-conditioned instances of fully real problems
can be done by homotopy methods.
The finding of 40 real solutions to the Stewart-Gough platform~\cite{die98}
is perhaps the most striking example.
The question is to find an efficient procedure to deform from the complex 
case to the fully real case.

\bigskip

\noindent {\bf Acknowledgments.}
The interest in homotopy methods by the FRISCO project has stimulated
the author's research.  The author is deeply indebted to all his co-authors.
The interactions with Birk Huber and Frank Sottile at MSRI were influential
in this current treatment of homotopies.

\newpage
{\small

}
  
\end{document}